\documentclass[12pt]{article}
\usepackage{amsfonts,amssymb,version}
\def\du{\coprod}

\def\P{\mbox{$\mathbb P$}}
\def\C{\mbox{$\mathbb C$}}

\def\E{{\cal E}}

\def\vee{\bigwedge}

\let\hra\hookrightarrow

\let\ra\rightarrow

\let\lra\longrightarrow

\let\ov\overline

\let \wt \widetilde

\newtheorem{theorem}{Theorem}[section]
\newtheorem{proposition}[theorem]{Proposition}
\newtheorem{lemma}[theorem]{Lemma}
\newtheorem{definition}[theorem]{Definition}
\newtheorem{corollary}[theorem]{Corollary}

\def\rem{\refstepcounter{theorem}\paragraph{Remark \thetheorem}}

\def\proof{\paragraph{Proof}}

\makeatletter \def\l@section{\@dottedtocline{1}{0em}{1.2em}} \makeatother

\begin{document}
\baselineskip=15pt

\title{Maximal subbundles and Gromov invariants}

\author{Yogish I. Holla}
\date{May 7, 2002}
\maketitle
\begin{abstract}
In this article we explicitly compute the number of maximal subbundles 
of rank $k$ of a generically stable bundle of rank $r$ and degree $d$ 
over a smooth projective curve $C$ of genus $g\ge 2$ over $\C$, 
when the dimension of the 
quot scheme of maximal subbundles is zero. Our method is to describe the 
this number purely in terms of the Gromov invariants of the 
Grassmannian and then use the formula of Vafa and Intriligator to compute them.
\end{abstract}

\section{Introduction}
Let $C$ be a smooth projective curve over $\C$
of genus $g\ge 2$.
For a vector bundle of rank $r$ and degree $d$ the $s$-invariant is defined by 
$s_k(E)=d \,k -r\,e_{{\rm max}}(E)$ where 
$e_{{\rm max}}(E)={\rm Max}\{{\rm deg}(F)\}$. 
Here the maximum is taken over all subbundles $F$ of $E$ of rank $k$.
It is known that for any $E$ the $s$-invariant satisfies 
$s_k(E)\le k(r-k)g$ (see the results of Mukai-Sakai \cite{Mukai-Sakai}).

One can get a more specific bound when we put additional structures on the 
bundle $E$. It is proved by Hirschowitz \cite{Hirschowitz} 
that for a general stable bundle 
the invariant $s_k(E)$ is independent of the choice of $E$ and satisfies 
 $s_{{\rm min},d}=s_k(E)=k(r-k)(g-1)+\epsilon$ where $\epsilon$ 
is the unique integer 
$0 \le \epsilon <r$ such that $s_{{\rm min},d}=kd ({\rm mod}\, r)$. 
Moreover if $e_{{\rm max},d}$ is the degree of a maximal subbundle of 
a general stable bundle $E$ of degree $d$ then every component of 
the quot scheme 
${\rm Quot}^{k,e_{{\rm max},d}}(E)$ of rank $k$ subsheaves of degree 
$e_{{\rm max},d}$ has dimension $\epsilon$.
And if the bundle is sufficiently general then the above quot scheme 
is itself smooth (see Lemma 2.1, Lange-Newstead \cite{Lange-Newstead} 
when the quot scheme is zero dimensional and in general see Remark 6.6, 
\cite{Holla}).  
Hence when $s_{{\rm min}}(d)=k(r-k)(g-1)$ we have a  zero dimensional 
quot scheme which is smooth. This defines the number $m(r,d,k,g)$ 
as the number of maximal subbundles of a sufficiently stable bundle $E$.

It is known that $m(r,d,1,g)=r^g$ 
(see Ghione \cite{Ghione}, Lange \cite{Lange3}, and Segre \cite{Segre} 
for $r=2$ and Okonek-Teleman \cite{Okonek-Teleman} and 
Oxbury \cite{Oxbury}, Theorem 3.1, in general). 
For the case when $k>1$ 
there is a formula to compute $m(r,d,k,g)$ when $k$ and 
$e_{{\rm max},d}$ are relatively prime (see Lange-Newstead 
\cite{Lange-Newstead}).

Our objective in this paper is to compute the number $m(r,d,k,g)$ 
explicitly for all choices of parameters (see Theorem \ref{formula}). 
We first define a notion of a twisted Gromov invariants for any vector bundle
following the methods of Bertram \cite{Bertram2}. The main tool this uses 
is the fact the quot schemes of our interested are generically smooth of 
expected dimension. This is the case when $s_e=d\,k-r\,e$ is sufficiently 
large. But if the bundle $E$ is sufficiently general then for all values of $e$
for which the quot scheme ${\rm Quot}^{k,e}(E)$ is non-empty it does 
satisfy the above property. This enable us to define the twisted Gromov 
invariants for all choices of $e \le e_{{\rm max}}(d)$. 
We then see that the numbers $m(r,d,k,g)$ can be computed purely in terms of 
these twisted Gromov invariants.
These invariants were also defined in \cite{Okonek-Teleman} and 
Behrend \cite{Behrend1}.

We show that these invariants are actually independent of the chosen 
general stable bundle. Also if $s_e$ is large enough then these invariants 
do not depend on the choice of the vector bundle of degree $d$. Then using 
Hecke transformations we compare the twisted Gromov invariants for 
different choices of $d$. These ideas will now enable us to compute 
the twisted Gromov invariants purely in terms of the Gromov invariants of 
the Grassmannian (when $E$ is trivial). And in this case there is a formula of 
Vafa and Intriligator 
(see \cite{vafa}, \cite{Intriligator}, and  for proof see Bertram 
\cite{Bertram2} and  \cite{Bertram4}) which exactly computes these invariants.
 
One observes that our formula for $m(r,d,k,g)$ and the one obtained in 
\cite{Lange-Newstead} are completely different in appearence. It will be 
interesting to compare the two expressions directly.

Some of these results can also be generalized to $G$-bundles for a 
connected reductive algebraic group $G$ and this will be done elsewhere. 

{\bf acknowledgment}  We wish  to thank G. Harder and M. S. 
Narasimhan who initiated him into this problem.
We also wish to thank Max-Planck Institute for Mathematics for hospitality.

\section{Twisted Gromov invariants}
Let $C$ be a smooth projective curve over $\C$. 
In this section we define the twisted Gromov invariants for a vector 
bundle $E$ of rank $r$ and degree $d$, and a positive integer $k\le r$.
Our view point here is that the construction of Gromov invariants given 
in Bertram \cite{Bertram2} holds for any vector bundle.

Let $E$ be a vector bundle of rank $r$ and degree $d$ over $C$.
We will be mostly be interested in the quot scheme
${\rm Quot}^{k,e}(E)$ parameterizing the rank $k$ subsheaves of degree $e$.
We have an open subscheme ${\rm Quot}^{k,e}_0(E)$ consisting of subbundles.
 We need some basic facts about these quot schemes which will be used in 
the sequel. We will denote by $s_e=d\,k-r\,e$.
First we consider the case $k=1$. 
Let $e_{{\rm max}}(E)$ be the degree of the maximal line subbundle of $E$.
We have natural morphisms  
$i_{e_1}:C_{(e_1-e)}\times {\rm Quot}^{1,e_1}(E) \ra {\rm Quot}^{1,e}(E)$ 
for each $e < e_1 \le e_{{\rm max}}$ which takes the pair
$(D,\, [L\hra E])$ to $[L(-D)\hra L\hra E]$. Using the universal properties of the quot schemes we check that this defines a morphism.
\begin{lemma}\label{emb}
The morphism $i_{e_1}$ is a closed embedding.
Moreover For $e\le e_1,\,e_2 \le e_{{\rm max}}$ with $e_1 \neq e_2$, we have
$$
i_{e_1}(C_{(e_1-e)}\times {\rm Quot}^{1,e_1}_0(E)) \cap 
i_{e_2}(C_{(e_2-e)}\times {\rm Quot}^{1,e_2}_0(E)) = \emptyset . 
$$ 
\end{lemma}
\proof If $[L\ra E]$ defines an element of 
${\rm Quot}^{1,e}(E)$ which is in image of $i_{e_1}$ then we can uniquely 
recover an effective divisor $D'$ and an element $[L'\hra E] \in 
{\rm Quot}^{1,e'}_0(E)$ for some $e'\ge e$ using the fact that 
$L'$ is the saturation of $L$ in $E$ and $D$ is  the divisor defined by 
the inclusion $L\hra L'$. This  fact already implies that the restriction 
of the morphism $i_{e_2}$ to 
$ C_{(e_2-e)} \times {\rm Quot}^{1,e_2}_0(E)$ is injective.
Now $i_{e_1}(D_1,\,[L_1 \hra E])=[L\hra E]$ if and only if 
the inclusion $L\hra L'$ factors via $L \hra L_1$ up to an isomorphism of 
$L$. Hence we see that 
the element defined by $L_1 \hra E$ in ${\rm Quot}^{1,e_1}_0(E)$ 
is in the image of corresponding morphism 
$C_{(e'-e_1)}\times {\rm Quot}^{1,e'}_0(E) \lra {\rm Quot}^{1,e_1}(E)$.  
This proves the second assertions of the lemma. $\hfill \square$

The above lemma defines a stratification on ${\rm Quot}^{1,e}(E)$
by locally closed subschemes $C_{(e_1-e)}\times {\rm Quot}^{1,e_1}(E)$, 
where $e\le e_1 \le e_{{\rm max}}$ 
and for $e_1=e$ we get an open strata ${\rm Quot}^{1,e}_0(E)$. 
Note that this stratification is in some sense weak as these quot schemes 
are not equidimensional and there can be stratas with $e_1 >e$ which 
might contain open subschemes of ${\rm Quot}^{1,e}(E)$.

Now we use the above lemma to get a similar understanding of the structure of 
${\rm Quot}^{k,e}(E)$ for $k>1$. For any subsheaf $F\hra E$ of rank 
$k$ we can take the $k$-th exterior power to get a line subsheaf 
$\vee ^k(F)\hra\vee ^k (E)$. This defines the Pl\"ucker morphism 
$P_e: {\rm Quot}^{k,e}(E)\ra {\rm Quot}^{1,e}(\vee ^k(E))$. 
Again it is easy to check that this morphism restricted to 
${\rm Quot}^{k,e}_0(E)$ is an embedding.
We have the following lemma which describes the boundary of 
${\rm Quot}^{k,e}(E)$ in terms of the stratification defined above.
\begin{lemma}\label{strat}
The scheme ${\rm Quot}^{k,e}(E)$ has a stratification 
${\rm Quot}^{k,e}_0(E)\du _{e\le e_1 \le e_{{\rm max}}} B_{e_1}$
where 
$B_{e_1} = P_e^{-1}(C_{(e_1-e)}\times {\rm Quot}^{1,e_1}_0(\vee ^k(E)))$. 
Moreover the image of the morphism $P_e$ when restricted to $B_{e_1}$ 
is exactly $C_{e_1-e}\times P_{e_1}({\rm Quot}^{k,e_1}_0(E)))$ 
and the fibers of 
$ {\rm pr}_2 \circ P_e$ when restricted to $B_{e_1}$ are irreducible and 
smooth of  dimension $k(e-e_1)$.
\end{lemma}
\proof The first assertion in the Lemma follows from the definition. 
The second follows from the fact that if $F\hra E$ is a rank $k$ 
subsheaf of $E$ then the saturation $F^{{\rm sat}}$ of $F$ in $E$ has the 
property that $\vee ^k(F^{{\rm sat}})$ is the saturation of $\vee ^k(F)$.
For the last part one observes that the fiber of the morphism 
${\rm pr}_2 \circ P_{e}|_{B_{e_1}}$ at the point  
$[F \hra E]$ is the quot scheme ${\rm Quot}^{k,e_1-e}(F)$ of 
length $e_1-e$ quotients ${\rm Quot}^{k,e_1-e}(F)$ and this space 
is irreducible and smooth of dimension $k(e_1-e)$. $\hfill \square$

\rem \label{est} The above lemma also gives us the following 
 dimension estimate
$$
{\rm dim}({\rm Quot}^{k,e_1}_0(E)) \le {\rm dim}({\rm Quot}^{k,e}(E))
-k(e_1-e).
$$ 
This estimate is also proved in Popa-Roth \cite{Popa-Roth}. 
It will be used to define the well definedness of the twisted Gromov 
invariants.

\rem \label{est1} If $Q$ is an irreducible  component of the quot scheme 
${\rm Quot}^{k,e}(E)$ such that every element in $Q$ 
corresponds to a quotient of $E$ which is not locally free then there 
exists an $e_1 \ge e$ such that $B_{e_1}$ is actually dense open subset of 
$Q$.This has the effect that we obtain a dimension bound
${\rm dim}(Q)\le {\rm Quot}^{k,e_1}_0(E) + k(e_1-e)$. 
This estimate will also be used later.

Let $x\in C$ be a point. We will denote by $E(x)$ the fiber of $E$ at $x$.
Let $Gr_k(E(x))$ denote the Grassmannian of $k$ dimensional subspaces of 
$E(x)$. Hence we have the universal  rank $k$ subbundle 
$ S(x) \ra E(x)\otimes {\cal O}_{Gr_k(E(x))}$.
There is a natural action of the general linear group $GL(E(x))$ of 
automorphisms of $E(x)$ on $Gr_k(E(x))$. 
Let $H\subset E(x)^*$ be a subspace of dimension $n\le k$. We then have a 
special Schubert variety
$Y_H$ defined by the degeneracy locus of the canonical map 
$H\otimes {\cal O}_{Gr_k(E(x))} \ra S(x)^*$. 
The variety $Y_H$ is irreducible and reduced of codimension 
$k+1-n$ and it represents the $(k+1-n)$-th Chern class of $S(x)^*$. 
For any element $g\in GL(E(x))$ 
we have a $g$ translate of $Y_H$ defined by the degeneracy locus associated 
to the translate $gH$.

For a Schubert variety $Y_H$ associated to $H$ we define the generalized 
Schubert scheme $W_e(x,H)=ev_x^{-1}(Y_H)$, where $ev_x$ is the evaluation map 
${\rm Quot}^{k,e}_0(E) \ra Gr_k(E(x))$ at $x$ which takes an element 
$[F\hra E]$ to the subspace $F(x)\hra E(x)$.

We need the following lemma.
\begin{lemma}\label{needed}
Given a vector bundle $E$ of rank $r$ over $C$ there is an integer
 ${\ov s}(E)$ such that for each $e$ with  
$s_e \ge {\ov s}(E)$, every component of the 
the quot scheme ${\rm Quot}^{k,e}(E)$ is generically smooth 
of expected dimension and a general element in every component 
corresponds to a subbundle of $E$.
\end{lemma}
\proof This is proved in  Popa-Roth \cite{Popa-Roth}, Theorem 5.14.  
In fact they prove a stronger version namely that there is an integer  
${\ov s}_1(E)$ such that for each $e$ with $s_e \ge {\ov s}_1(E)$ the
quot scheme ${\rm Quot}^{k,e}(E)$ is also irreducible. $\hfill \square$

The above lemma allows us to define the invariant ${\ov s}(E)$ 
(namely minimal such integer satisfying the above lemma) for any 
vector bundle.

Now we define the twisted Gromov invariants as the 
intersection numbers in the quot scheme. 
Let $e$ be such that $s_e \ge {\ov s}(E)$.
Let $(x,H)$ be a pair with $x \in C$ and $H$   
a subspace of $E(x)^*$ of rank $n \le k$. 
We define the subscheme  $V_e(x,\,H) \subset 
{\rm Quot}^{k,e}(E) $
as the the degeneracy locus of the composite morphism 
$H\otimes {\cal O}_{{\rm Quot}^{k,e}(E)} \ra ({\rm pr}_1 ^*E)(x) 
\ra {\cal F}^*(x)$, where
$ ({\cal F}^*)(x)$ is the restriction of the dual of the universal subsheaf
${\cal F}\subset {\rm pr}_1 ^*E$ of rank $k$ and degree $e$ to 
$\{x\}\times {\rm Quot}^{k,e}(E)$ and 
${\rm pr}_1: C\times {\rm Quot}^{k,e}(E) \ra C$ is the first projection.
By definitions we can check that $W_e(x,\, H)$ is the scheme theoretic 
intersection of $V_e(x,\,H)$ and ${\rm Quot}^{k,e}_0(E)$.
\begin{lemma}\label{import}
With the notations above if $g\in GL_{x}(E)$ general then 
every component of $V_e(x,\,gH)$ has codimension exactly $k-n+1$ in 
${\rm Quot}^{k,e}(E)$ and $V_e(x,\,gH)$ represents the 
$(k-n+1)$-th Chern class of ${\cal F}^*(x)$.
\end{lemma}
\proof
Firstly the conditions $s_e \ge {\ov s}(E)$ ensures that the quot scheme 
${\rm Quot}^{k,e}(E)$ has the properties mentioned in the Lemma \ref{needed}
This has the effect that the above quot scheme is locally of complete 
intersection hence Cohen-Macaulay 
(see for example K\"ollar, I, Theorem 5.17, \cite{Kollar}). 
Hence the subscheme $V_e(x,\,H)$ represents the 
$k-n+1$-th Chern class of ${\cal F}^*(x)$ if we show that the codimension 
of every irreducible component of $V_e(x,\,H)$ is $k-n+1$ 
(see Fulton,  Theorem 14.3, Fulton \cite{Fulton}). 
Now we use the Lemma \ref{strat} \and Lemma \ref{needed} 
we can continue as in the proof of 
the first part of the Theorem 1.4 of  \cite{Bertram2} to check the
codimensionality condition. $\hfill \square$

The above lemma allows us to define the twisted Gromov invariants as follows.
Let $E$ be a rank $r$ and degree $d$ vector bundle over $C$.
Let ${\ov s}(E)$ be as defined in Lemma \ref{needed}.
Let $X_1,\,\ldots,\, X_k$ be weighted variables such that weight of $X_i$ is 
$i$. Let $P(X_1,\,\ldots,\, X_k)$ be a homogeneous polynomial of weighted 
degree $s_e +k(r-k)(1-g)$ with $s_e=d\,k-r\,e$ and $s_e \ge {\ov s}(E)$.
\begin{definition} We define the
twisted Gromov invariant  
$N_{d,e}(P(X_1,\,\ldots,\, X_k),\,E)$ for a vector bundle $E$, and an integer 
$e$  such that $s_e > {\ov s}(E)$,
as the intersection number 
$$
P(c_1({\cal F}^*(x),\,\ldots ,\,c_k({\cal F}^*(x))[{\rm Quot}^{k,e}(E)]
$$
\end{definition}
Here $[{\rm Quot}^{k,e}(E)]$ denotes the fundamental cycle.
\rem The above number is independent of the single point $x$ chosen. More 
generally we could choose points $x_i$ for $i=1,\,\ldots,\,n$ and subspaces 
$H_i \subset E(x_i)^*$ of dimension $a_i\le k$ such that 
$\sum _i (k-a_i+1) =s_e +k(r-k)(1-g)$ then 
the number $N_{d,e}(\Pi_{i=1}^n X_{k-a_i+1},\,E)$ can be exactly computed 
as the intersection number
$\cap V_e(x_i,\,g_i H_i)$ for general elements $g_i\in GL(E(x_i))$.

The following proposition shows how to compute the above intersection as 
an intersection in the open subscheme ${\rm Quot}^{k,e}_0(E)$.

\begin{proposition}\label{welldefined}Let $E$ be a rank $r$ vector bundle of 
degree $d$ and let $e$ be such that $s_e \ge {\ov s}(E)$. 
For $i=1,\,\ldots ,\,N$, let 
$x_i \in C$ be distinct points and $H_i \subset E^*(x_i)$ be subspaces 
of dimension $a_i\le k$ such that $\sum _{i=1}^N (k-a_i+1) =s_e +k(r-k)(1-g)$. 
Then for general choices of $g_i\in GL(E(x_i))$, the twisted Gromov invariants 
$N_{d,e}(\Pi_{i=1}^N X_{k-a_i+1},\,E)$ can be exactly computed 
as the number of intersections $\cap _{i=1}^NW_e(x_i,\,g_iH_i)$ 
(counted with multiplicities).
\end{proposition}
\proof 
The proof of the proposition exactly follows second part of the proof 
of Theorem 1.4 of  \cite{Bertram2}. Here we can
use the dimension bound obtained in the Remark \ref{est} to ensure that 
for general translates the intersections in the boundary are trivial for 
any choice of $e$ such that $s_e \ge {\ov s}(E)$. $\hfill \square$

Now we compare the intersection numbers obtained in the above theorem  
for different rank $r$ vector bundles of a fixed degree $d$.

Let $E_1$ and $E_2$ be two rank $r$ vector bundles over $C$ of degree $d$.
It can be proved using the smoothness and irreducibility of the moduli 
stack bundles of rank $r$ and degree $d$ that there is a smooth irreducible 
variety $B$, a family
of bundles $\E$ over $C\times B$ and two points $x_1$ and $x_2$ of $B$ 
such that $\E |_{C\times x_i}\cong E_i$ for $i=1,\,2$. 
We will denote by $\E _x$ 
the bundle $\E |_{C\times x}$ for $x\in B$. 
We need the following generalization of the Lemma \ref{needed}.
\begin{lemma}\label{irred}
There exists an integer $s_{{\cal E}}$ independent of $x\in B$ 
such that for all $s\ge {\ov s}_{{\cal E}}$ the quot scheme 
${\rm Quot}^{k,e}(\E_x)$ is generically smooth of 
expected dimension and satisfies the property that general elements in
every irreducible component of ${\rm Quot}^{k,e}(\E_x)$
lies in ${\rm Quot}^{k,e}_0(\E_x)$
\end{lemma}
\proof The lemma follows from the arguments similar to the proof of 
Theorem 6.4 of \cite{Popa-Roth} once we establish the existence of an 
integer ${\ov s}_{{\cal E}}$ independent of $x\in B$ such that 
${\rm Quot}^{k,e}_0(\E_x)$ is generically smooth of expected dimension 
for all $e$ with  $s_e \ge {\ov s}_{{\cal E}}$.
The last statement can be again proved by going through the arguments 
of Proposition 6.1 and Theorem 6.2 of \cite{Popa-Roth} 
(Also see Proposition 5.11, \cite{Holla} where the integer 
${\ov s}_{{\cal E}}$ is described 
in terms of instability degree of a principal bundle). $\hfill \square$

We can use the above Lemma to define the invariant ${\ov s}_{{\cal E}}$
(namely minimal such integer satisfying the above lemma) 
for a family of vector bundles ${\cal E}$ over $C\times B$.

\begin{proposition}\label{well} If $\E$ is a family of rank $r$ vector 
bundles of degree $d$ on 
$C\times B_1$ with $B_1$ being a  smooth curve and $e$ is chosen such that 
$s_e \ge {\ov s}_{{\cal E}}$ then 
the Gromov invariants $N_{d,e}(P(X_1,\,\ldots,\,X_k),\E _x)$ 
are independent of the choice of points $x\in B$.
\end{proposition}
\proof We consider the relative Quot scheme 
$f:{\rm Quot}^{k,e}(\E)\ra B_1$ which has the property that for each 
$x\in B_1$ the fiber is exactly ${\rm Quot}^{k,e}(\E _x)$. If
$s_e \ge {\ov s}_{{\cal E}}$ then
by Lemma \ref{irred} we see that each of ${\rm Quot}^{k,e}(\E _x)$ are 
generically smooth of expected dimension. Hence 
$f$ is a locally complete intersection morphism 
(see I, Theorem 5.17, \cite{Kollar}), 
in particular flat. Now the proposition follows from Lemma 1.6 of 
\cite{Bertram2}. $\hfill \square$ 

The above proof can also be used to show that the intersection numbers are 
actually depend on only the genus $g$ of the smooth projective curve $C$
(for the case of trivial bundles, this is the Proposition 1.5 \cite{Bertram2}).

We also have a simple relation between the twisted Gromov invariants when the 
vector bundle is tensored with a line bundle.

\begin{lemma}\label{tensor}
Let $E$ be a rank $r$ vector bundle of degree $d$ and $L$, a line bundle 
of degree $d_1$. Then we have ${\ov s}(E)={\ov s}(E\otimes L)$ and the
twisted Gromov invariants of $E$ and $E\otimes L$ are related by the 
following
$$
N_{d+rd_1,e+kd_1}(P(X_1,\,\ldots,\,X_k),\,E\otimes L)=
N_{d,e}(P(X_1,\,\ldots,\,X_k),E).
$$ 
\end{lemma}
\proof This follows from the fact that the quot schemes 
${\rm Quot}^{k,e+kd_1}(E\otimes L)$ and ${\rm Quot}^{k,e}(E)$ are naturally 
isomorphic and the isomorphism preserves the all the degeneracy loci.
$\hfill \square$

\section{Generically stable bundles}
In this section we will assume that $C$ is a smooth projective curve 
over $\C$ of genus $g \ge 2$.

we first prove some basic facts about the structure of the quot schemes for 
generically stable bundles and these fact will be used later.

We first begin by the following result which will enable us to define 
twisted Gromov invariants when the degrees of the subbundles is not very small.
Let $M^s(r,d)$ denote the coarse moduli space of stable bundles of rank $r$ 
and degree $d$. Let $0< k <r$. 
For a vector bundle $E$ over $C$ of rank $r$ and degree $d$, 
recall the notion of the $s$ invariant by 
$s_k(E)=d \,k -r\,e_{\rm max}(E)$ where 
$e_{\rm max}={\rm Max}\{{\rm deg}(F)\}$.
Here the maximum is taken over all rank $k$ subbundles $F$ of $E$.
It is known that if the vector bundle $E$ is generically stable then 
$s_{{\rm min},d}=s_k(E)=k(r-k)(g-1)+\epsilon$ where $\epsilon$ is the
unique integer with $0\le \epsilon \le r-1$ such that $s=k\,d ({\rm mod}\, r)$.
Let $e_{{\rm max},d}$ be the degree of the maximal subbundle of a 
generically stable bundle of degree $d$.

\begin{proposition} \label{generic}
There exists a non empty open subset $U({r,d})\subset M(r,d)$ with the 
property that for each $E\in U({r,d})$ and for every $e\le e_{{\rm max},d}$, 
every component of quot scheme 
${\rm Quot}^{k,e}(E)$ is generically smooth of expected 
dimension $(=d\,k-r\,e+k(r-1)(g-1))$ and satisfies the property that 
general elements in every irreducible component corresponds to 
subbundles of $E$.
\end{proposition} 
\proof
It is already known that there is a non empty open subset $U(r,d)$ of
$M(r,d)$ such that each $E$ in $U(r,d)$ satisfies 
$s_{{\rm min},d}=s_k(E)=d\,k-r\,e_{{\rm max},d}$ and the dimension estimate 
${\rm dim}({\rm Quot}^{k,e}(E))\le s_e +k(r-k)(1-g)$ holds for all 
$e \le e_{{\rm max}}$ (see Example 5.4 \cite{Popa-Roth}).
This along with the deformation theoretic lower bounds implies that every 
component of the quot scheme has the above dimension.
Moreover it is also known that every irreducible component of 
${\rm Quot}^{k,e}_0(E)$ is generically smooth (see for example 
Proposition 6.8, \cite{Holla}).
So we only have to only
make sure that these quot schemes have no pathological components.
Let $Q$ be an irreducible component of ${\rm Quot}^{k,e}(E)$ such that 
every element in $Q$ corresponds to quotient which is not torsion free. 
Then by Remark \ref{est1} there is an $e_1 > e$ such that 
${\rm dim}(Q)\le {\rm dim}({\rm Quot}^{k,e_1}_0(E)) + k(e_1-e)$. Now this 
leads to a contradiction once we put the values of the dimensions of 
these spaces.   $\hfill \square$

The above result shows that if $E\in U(r,d)$ then 
${\ov s}(E)=s_{{\rm min},d}$. This enables us to define the twisted
Gromov invariant $N_{d,e}(P(X_1,\,\ldots,\,X_k),E)$ for all possible values 
of $e \le e_{{\rm max},d}$.
Moreover if we stay in the open subset $U(r,d)$ then 
we see by Proposition \ref{well} that these numbers are actually independent 
of the choice of $E\in U(r,d)$. This enable us to define 
$N_{d,e}(P(X_1,\,\ldots,\,X_n))$ to be $N_{d,e}(P(X_1,\,\ldots,\,X_n),E)$ 
for some $E \in U(r,d)$ and $e$ with $e\le e_{{\rm max},d}$. 

Also the Proposition  \ref{well} shows that if $E_1$ is any vector bundle 
of rank $r$ and degree $d$ then for any $e$ with $s_e \ge {\ov s}(E_1)$ 
we have 
$$
N_{d,e}(P(X_1,\,\ldots,\,X_n))=N_{d,e}(P(X_1,\,\ldots,\,X_n),E_1).
$$

Now we will compare the above invariants when the degree of the vector 
bundle changes. The main result is the following.

\begin{proposition}\label{compare}
For a polynomial $P(X_1,\,\ldots,\,X_k)$ of 
weighted degree $s_e +k(r-k)(g-1)-k$ with $e$ satisfying 
$e\le e_{{\rm max},d-1}$ We have 
$$
N_{d-1,e}(P(X_1,\,\ldots,\,X_k))=N_{d,e}(X_kP(X_1,\,\ldots,\,X_k))
$$
\end{proposition}
We first recall a basic lemma which will be a step in the proof of the 
above result.
Let $x\in C$ be a point fixed. For a bundle $E$ of rank $r$ and degree $d$
and a quotient space  $l:E(x)\ra \C$ (or equivalently a line 
$l\subset E^*(x)$) of dimension $1$ at a point $x$ we  
have the Hecke transform $F_l ={\rm ker}({\wt l})$, where 
${\wt l}$ is the map $E\ra \C$ defined by $l$. Then the degree of $F_l$ is 
equal to $d-1$. Hence for any point $g\in GL(E(x))$ 
we have a vector bundle $F_{gl}$ defined by the above procedure. 

\begin{lemma}\label{hecke}
For any two non-empty open sets $U \subset M^s(r,d)$ and 
$U_1\subset M^s(r,d-1)$   
there exists a vector bundle $E\in U$ and a non empty open set 
$V\subset GL(E(x))$ such that for each $g\in V$ the 
vector bundle  $F_{gl}$ got by Hecke transform with respect to 
$gl$ lies in $U_1$.
\end{lemma}
\proof 
Let ${\rm Vect}(r,d)$ denote the moduli stack of vector bundles of rank $r$ 
and degree $d$.
We have the Hecke stack ${\cal H}$ defined by triples $(F,E,i)$ where 
$F$ and $E$ are rank $r$ bundles of degree $d-1$ and $d$ respectively 
and an inclusion $i:F \le E$ whose cokernel is supported at a single point $x$.
We also have two morphisms $h_1:{\cal H} \ra {\rm Vect}(r,d-1)$ and 
 $h_2:{\cal H} \ra {\rm Vect}(r,d)$ which takes such a triple to $F$ and 
$E$ respectively. It is easy to check that these two morphisms are 
surjective. Now the stack $H$ is irreducible as  
is a projective bundle over ${\rm Vect}(r,d)$, which is irreducible, 
whose fiber at $E$ is given by $Gr_m(E(x))$. 
Since  $M^s(r,d)$ and  $M^s(r,d-1)$ are the coarse moduli spaces of 
open substack of ${\rm Vect}(r,d)$ and ${\rm Vect}(r,d-1)$ respectively,
hence we can check that the inverse images of the open sets $U$ and $U_1$
under the morphism $h_1$ and $h_2$ respectively intersect in a non empty 
open substack of ${\cal H}$. This proves the lemma. $\hfill \square$

Now we return to the proof of the Proposition \ref{compare}. 
Let $x\in C$ be a fixed point.
Choose a vector bundle $E \in U(r,d)$ prescribed by the 
Lemma \ref{hecke} for the open subsets $U_1=U(r,d-1)$ and  
$U=U(r,d)$. We then obtain a non empty open subset
$V \subset GL(E(x))$ such that for each $g\in V$ the bundle $F_{gl}$ lies in 
$U_1$. 
Let $x_1,\,\ldots,x_N$ be distinct points of the curve $C$ 
which are also distinct from $x$ and $H_i \subset E^*(X_i)$
be subspaces of rank $a_i \le k$ such that 
$\sum _{i=1}^N (k-a_i+1) =s_e +k(r-k)(1-g)-k$.
Then by Proposition \ref{well} the twisted Gromov invariant 
$N_{d,e}(X_k\Pi_{i=1}^l X_{k-a_i+1},\,E)$ 
can be exactly computed 
as the number of intersections 
$V(x,\,gl)\bigcap _{i=1}^NV_e(x_i,\,g_iH_i)$ for general elements of $g_i$'s 
and $g$, where  $V_e(x_i,\,g_iH_i)$ (and $V(x,\,gl)$) are the degeneracy loci
in ${\rm Quot}^{k,e}(E)$ associated to $H_i$ (and $l$). 

We have a natural morphism 
$f:{\rm Quot}^{k,e}(F_{gl})\ra {\rm Quot}^{k,e}(E) $ which takes 
$[H\hra F_{gl}]$ to $[H\hra F _{gl}\hra E]$. This morphism can be  
checked to be an embedding.

Now we check by definitions that the degeneracy locus $V(x,\,gl)$ is 
exactly the isomorphic image of $f$ and for each $i$ the scheme  
$f^{-1}V_e(x_i,\,g_iH_i)$  is exactly the degeneracy locus associated to 
$F_{gl}$ for the subspace $gH_i \subset F_{gl}^*(x_i)$ 
Hence for general $g_i$ and $g$ there is a natural isomorphism between the 
the zero dimensional schemes 
$$
V(x,\,gl)\cap _{i=1}^NV_e(x_i,\,g_iH_i) \cong 
\cap _{i=1}^Nf^{-1}(V_e(x_i,\,g_iH_i)).
$$ 
Hence by Proposition \ref{welldefined} we see that the left hand side of 
the above isomorphism computes the twisted Gromov invariants for $E$ and 
the right hand side of the of the above isomorphism computes the twisted 
Gromov invariants for $F_{gl}$.  
This completes the proof of the Proposition \ref{compare}. 
$\hfill \square$       

\rem the last part of the proof of the above proposition also follows from 
the excess intersection formula (see Fulton, Theorem 6.3 and Proposition 6.6,
\cite{Fulton}) since 
the quot schemes in question satisfies the local complete intersection
property. 

Recall the definition of the standard Gromov invariants 
$N_{e}(P(X_1,\,\ldots,\,X_k),\,g)$ for the Grassmannian
from  \cite{Bertram2} for a polynomial 
$P(X_1,\,\ldots,\,X_k)$ of weighted degree $e\,r+k(r-k)(g-1)$ with $-e$ 
large enough.
This invariant in our notations is equal to 
$N_{0,-e}(P(X_1,\,\ldots,\,X_k))$ for the choices of 
$e$.
It is also proved in Lemma 5.3, \cite{Bertram} 
(Also see \cite{Bertram2}, Remark before the Proposition 1.7) that the 
above invariant can be 
consistantly defined for larger values of $e$ by the following
is defined for larger $e$ by the recurssion 
$N_{0,e}(P(X_1,\,\ldots,\,X_k))=N_{0,e-k}(X_k^rP(X_1,\,\ldots,\,X_k))$.
This equality also follows from Lemma \ref{tensor}.

We have the following explicit formula for the twisted Gromov 
invariants in terms of the standard Gromov invariants of the Grassmannian.

\begin{theorem}\label{main}
Let $r$ and $k$ be fixed. Let $d=ar-b$ with $0\le b <r $
and $e\le e_{{\rm max}}(d)$. Let $P(X_1,\,\ldots,\,X_k)$ be a polynomial 
of weighted degree $d\,k-r\,e+k(r-k)(1-g)$. 
Then all the twisted Gromov invariants 
are computable in terms of the actual Gromov invariants of the Grassmannian 
by the following 
$$
N_{d,e}(P(X_1,\,\ldots,\,X_k))=N_{0,e-ak}(X_k^{b}P(X_1,\,\ldots,\,X_k))
$$
\end{theorem}
\proof This just follows from the 
Proposition \ref{compare} and Lemma \ref{tensor}.  $\hfill \square$

\rem The above result shows that the twisted Gromov invariants are independent 
of the choice of the genus $g$ curve $C$.

Now we record here the formula of Vafa and Intriligator, proved by 
A. Bertram (see \cite{Bertram2}, \cite{Bertram4}) and Siebert-Tian 
(see \cite{Siebert-Tian}) about explicit 
computation of Gromov invariants $N_{0,e}(P(X_1,\,\ldots,\,X_k))$.  

Let $P(X_1,\,\ldots,\,X_k)=\Pi_{i=1}^mX_{a(i)}$ be a polynomial with 
$0<a_i\le k$ such that the weighted degree of $P$ is 
$\sum_i (k-a_i+1)=-er+k(r-k)(1-g)$. Then we have the following.
\begin{proposition}\label{vafa} For the polynomial $P=\Pi_{i=1}^mX_{a(i)}$
as above, 
the Gromov invariant 
$N_{0,e}(P(X_1,\,\ldots,\,X_k))$ is given by the following.
$$ 
\frac{r^{k(g-1)}(-1)^{-e(k-1)+(g-1)k(k-1)/2}}{k!}
\sum_{\{(\rho_1,\,\ldots,\,\rho_k)|\rho _i^r=1;\rho _i\neq \rho _j\}}
\frac{\left(\Pi_{l=1}^m \sigma_{k-a(l)+1}(\rho)\right)}
{\left(\Pi_{i=1}^r{\rho_i}\,\Pi _{i\neq j}(\rho_i-\rho_j)\right)^{(g-1)}}
$$
where $\sigma_{j}(\rho)$ is the $j$-th symmetric polynomial in $\rho_i$'s. 
\end{proposition}
 
The above proposition along with the Theorem \ref{main} we can  
explicitly compute the twisted Gromov invariants.

\section{Maximal Subbundles}

In this section we will relate the twisted Gromov invariants to the number of 
Maximal subbundles of a generically stable bundles.
We will assume that the genus of the curve is atleast $2$.
Recall that for a generically stable bundle $E$ of rank $r$ and degree 
$d$ we have 
$s_{{\rm min},d}=s_k(E)=k(r-k)(g-1)+\epsilon$ where $\epsilon$ is a 
unique integer 
$0\le \epsilon \le r-1$ such that $s_{{\rm min},d}=kd ({\rm mod}\, r)$.
Also recall the definition of $e_{{\rm max},d}$ from the previous section.

We need the following proposition in the sequel. 

\begin{proposition}\label{smooth} Let $r$, $k$, $d$ and $e_{{\rm max},d}$ 
be as before. Then for a sufficiently general stable bundle $E$ the quot 
scheme ${\rm Quot}^{k,e_{{\rm max},d}}(E)$ is a smooth scheme of dimension 
$s_{{\rm min},d}+k(r-k)(1-g)$ consisting only of vector bundle quotients  
\end{proposition}
\proof We have already seen that for generically stable bundle 
$E$ every component of the quot 
scheme ${\rm Quot}^{k,e_{{\rm max}}}(E)$ is of expected dimension. 
Now the result of Mukai and Sakai \cite{Mukai-Sakai}, as $e=e_{{\rm max}}$, 
ensures there is no boundary. Hence we have to only show that 
that for a sufficiently general $E$ the  scheme 
 ${\rm Quot}^{k,e_{{\rm max}}}_0(E)= {\rm Quot}^{k,e_{{\rm max}}}(E)$ 
is smooth. This follows from Remark 5.4, of Holla \cite{Holla} where this 
has been worked out for a principal $G$-bundle with $G$ a connected 
reductive algebraic group. $\hfill \square$

\rem The above result follows from Lemma 2.1 of \cite{Lange-Newstead} 
for the case when  $s_{{\rm min},d}=k(r-k)(g-1)$. 
This is the only case when we will use the Proposition \ref{smooth} 
Also the above result has been worked out in \cite{Oxbury}, Proposition 1.4 
when $k=1$. 

\rem The Theorem \ref{main} allows us to compute explicitly 
all the Chern numbers of the vector bundle ${\cal F}^*(x)$, where 
${\cal F}$ is the universal subbundle of ${\rm pr}_1^*E$ on 
$C\times {\rm Quot}^{k,e_{{\rm max},d}}(E)$ interms of the 
Gromov invariants. From here one should be able 
to compute all the Chern numbers of the quot scheme itself. 
Since the Proposition \ref{smooth} ensures that these quot schemes are 
smooth projective schemes hence it would be interesting to know what 
kind of objects these are. For the case $k=1$ it is a conjecture of 
Oxbury (see \cite{Oxbury}, Conjecture 2.8) that these quot schemes are 
irreducible when they are of dimension greater than zero. 
We ask whether this conjecture is true for other $k$.

Now we will handle the case when the choice of $k$, $r$ and $d$ are made 
such that $k(r-k)(g-1)=k\,d ({\rm mod}\, r)$. In this case we can write 
down $e_{{\rm max},d}=dk+k(r-k)(1-g)$. This has the effect that 
for a bundle $E$ which is generically stable in the sense defined in 
Proposition \ref{smooth}, then the quot scheme 
${\rm Quot}^{k,e_{{\rm max},d}}(E)$ is of dimension $0$ and smooth. 
  
Hence we can count the number points in the quot scheme. 
We denote this number by $m(r,d,k,g)$.

We have the following explicit formula for the number $m(r,d,k,g)$.

\begin{theorem}\label{formula}Let $d=ar-b$ with $0\le b <r$.
The number $m(r,d,k,g)$ is calculated by the 
following formula.
$$ 
\frac{r^{k(g-1)}(-1)^{(k-1)(bk-(g-1)k^2)/r}}{k!}
\sum_{(\rho_1,\,\ldots,\,\rho_k)|\rho _i^r=1;\rho _i\neq \rho _j}
\frac{\left(\Pi_{l=1}^m \rho_l \right)^{b-g+1}}
{\left(\Pi _{i\neq j}(\rho_i-\rho_j)\right)^{(g-1)}}.
$$
\end{theorem}
\proof This follows from the fact that the number $m(r,d,k,g)$ we are 
interested in is exactly $N_{d,e_{{\rm max}},d}(1)$ for our choice 
of vector bundle $E$ in Proposition \ref{smooth} as the expected dimension 
of ${\rm Quot}^{k,e_{{\rm max}}}(E)$ is $0$. 
Now using the Theorem \ref{main}
we obtain $m(r,d,k,g)=N_{0,e_{{\rm max}}-ak}(X_k^{b})$.
Now the Proposition \ref{vafa} implies the theorem.  
$\hfill \square$

\begin{corollary}
$m(r,d,1,g)=r^g$
\end{corollary}

The above result recovers a part of the Theorem 3.1 of Oxbury \cite{Oxbury}.

The next case is when $k=2$.
We write $d=ar-b$ and in this case the zero dimensionality of the quot 
scheme gives us the condition $2b-4(g-1)=0 ({\rm mod}\, r)$. Now using the 
fact that $\sum _{\{\rho ^r=1,\rho\neq 1\}}\rho^i$ is equal to $-1$, 
if $0<i<r$ and is equal to $r-1$ if $i=0,r$,  
the Theorem \ref{formula} gives 
$$ m(r,d,2,g)=
\frac{r^{2(g-1)+1}(-1)^{g-1+(2b-4(g-1))/r}}{2}
\sum_{\{z|z^r=1,z\neq 1\}}
\frac{z^{b-g+1}}{(1-z)^{2(g-1)}}.
$$
The last summation can be computed using some combinatorial
identities which we will now briefly describe.
We define $B(k,m)=\sum_{\{z|z^r=1,z\neq 0\}}z^m(1-z)^{-k}$. 
Then we can compute 
$B(k,m)$ using the recursive relation $B(k,m)=B(k,m-1)+B(k-1,m-1)$
and the condition $\sum_{i=0}^{r-1}B(k,i)=0$.
These relations can easily be verified and one can hence obtain the following
recursive formula for $m>0$ 
$$
B(k,m)= B(k,0)-\sum_{i=0}^{m-1}B(k-1,i)
$$
and $B(k,0)$ can be computed by the following formula.
$$
B(k,0)=\frac{1}{r}\sum_{i=0}^{r-2}(r-i-1)B(k-1,i)
$$
For example using the above formulas we can calculate
for $0<m<r-1$
$$
B(2,m)=-(r^2+r(6-6m)+6m^2-12m+5)/12.
$$
Let $r\ge 3$.
If we assume that the genus of the curve is $2$, then the condition 
for the quot scheme to be zero dimensional is $2b-4=0 ({\rm mod}\, r)$. 
Since $0\le b <r$, hence either $b=2$ or 
$2b-4=r$. Hence we have the following corollary
\begin{corollary}Let $g=2$ and $r\ge 3$. 
If $b=2$ then we have $m(r,d,2,2)=r^3(r^2-1)/24$
and if $r=2b-4$ then we have $m(r,d,2,2)=r^3(r^2+2)/48 $.
\end{corollary}
The formula for $r=2b-4$ was also obtained by Lange and Newstead 
(Theorem 4.1 \cite{Lange-Newstead}).


\end{document}